\documentclass{amsart}
\usepackage{graphicx}
\usepackage{latexsym}
\usepackage{amsmath}
\usepackage{amssymb}
\usepackage{amsthm}
\usepackage{amscd}
\usepackage{amsrefs}

\newtheorem{thm}{Theorem}[section]

\newtheorem{prop}[thm]{Proposition}
\newtheorem{cor}[thm]{Corollary}

\newtheorem*{thms}{Theorem}

\theoremstyle{definition}
\newtheorem{definition}[thm]{Definition}
\newtheorem{remark}[thm]{Remark}

\newtheorem{app}[thm]{Application}

\newcommand{\al}{\alpha}                
\newcommand{\sig}{\sigma}               

\newcommand{\ra}{\rightarrow}           
\newcommand{\lra}{\longrightarrow}

\newcommand{\subs}{\subseteq}           

\newcommand{\sect}{{\mathcal{x}}}
\newcommand{\e}{{\acute{\textup{e}}}}
\newcommand{\et}{{\acute{e}t}}
\newcommand{\Z}{{\mathbb{Z}}}
\newcommand{\Ze}{{\Z_{\et}}}
\newcommand{\R}{{\textup{R}}}

\newcommand{\Q}{{\mathbb{Q}}}
\newcommand{\G}{{\mathbb{G}}}

\newcommand{\Hc}{{\mathcal{H}}}
\newcommand{\HT}{{\widehat{H}}}

\newcommand{\trg}{{\tau_{_{\leq 0}} \textup{R}\Gamma}}
\newcommand{\RG}{{\textup{R}\Gamma}}
\newcommand{\tn}{{\tau_{_{\leq 0}}}}

\newcommand{\Gal}{{\textup{Gal}}}
\newcommand{\Hom}{{\textup{Hom}}}

\newcommand{\Spec}{{\textup{Spec} \ }}

\newcommand{\gr}{{\textup{Gr}}}
\newcommand{\K}{{\mathcal{K}}}

\author{Karen Acquista}
\title{A criterion for cohomological dimension}
\date{\today}
\email{kea@math.bu.edu}
\address{Boston University \\ Department of Mathematics and Statistics \\ Boston, MA}

\begin{document}

\begin{abstract}  We give a criterion for the cohomological dimension of a field, involving norm maps on Milnor $K$-theory; this criterion was 
originally formulated by Kato.  The theorem we prove is a generalization 
of a theorem in Serre's book on Galois cohomology. 
\end{abstract}

\maketitle

Let $G$ be a profinite group.  $G$ is said to have cohomological dimension at most $n$ if for every discrete torsion $G$-module $A$, 
$$
H^q(G, A) = 0
$$
for all $q>n$, see $\sect$I.3 of Serre \cite{CG}; the cohomological dimension of $G$ will be denoted $cd(G)$.    

For a field $k$, one sets the cohomological dimension of $k$ to be the cohomological dimension of the absolute Galois group $G_k$.  There is a classical criterion for fields of cohomological dimension of most one:

\begin{thms}[Serre  \cite{CG}*{$\sect$II.3}] Let $k$ be a perfect  field.  Then the cohomological dimension of $k$ is at most 1 if and only if for every finite extension $K/k$ and every finite Galois extension $L/K$, the norm map:
$$
N : L^{\times} \lra K^{\times}
$$
is surjective.  
\end{thms}

The purpose of this note is to use motivic complexes to obtain an analogous criterion for fields of cohomological dimension at most $n$.    

\section{The main result}

For a field $K$, let $K_n(K)$ denote the $n^{th}$ Milnor $K$-group of $K$. It is defined as:
$$
K_n(K) = K^{\times} \otimes \ldots \otimes K^{\times}/I,
$$
where $I$ is the ideal generated by $a_1 \otimes \ldots \otimes a_n$ such that $a_i + a_j = 1$ for some $i \neq j$.  

Let $L/K$ be a finite Galois extension.  Although $K_n(L)^{\Gal(L/K)}$ is not in general isomorphic to $K_n(K)$ for $n \geq 2$, one can still define a norm map:
$$
N : K_n(L) \lra K_n(K),
$$
see Bass and Tate \cite{BT}.  It was shown by Kato \cite{K2}*{$\sect$1} 
that these norm maps satisfy certain functorial properties which 
characteruze them uniquely, see $\sect$IX.3 of Fesenko-Vostokov \cite{FV} 
for an overview.  One can use motivic complexes to construct the norm map simply using the Galois action, see Remark \ref{norms}.  

For the convenience of the reader, the properties of motivic complexes used in the course of the proof of the main result are included in Appendix \ref{motivic}.  We will also need some facts about Tate cohomology with coefficients in a complex; these can be found in Appendix \ref{tate}.

We have the following result:

\begin{thm} Let $k$ be a field; if the characteristic of $k$ is $p>0$, assume furthermore that $[k:k^p] \leq p^{n-1}$.  Then the cohomological dimension of $k$ is at most $n$ if and only if for every finite extension $K/k$ and every finite Galois extension $L/K$, the norm map:
$$
N : K_n(L) \lra K_n(K)
$$
is surjective.
\end{thm}

This criterion for $k$ was originally formulated by Kato \cite{K2}; in $\sect$3.3 Kato calls such a field $k$ a $\mathfrak{B}_n$-field.  

\begin{proof}  Let $k$ be a field such that for all finite extensions $K/k$ and all finite Galois extensions $L/K$, the norm map:
$$
N : K_n(L) \lra K_n(K)
$$
is surjective.  

Let $\Z_B(n, L)$ denote the weight $n$ motivic complex for $L$ over the Zariski site, see Remark \ref{zb}.  Our assumption implies that for all intermediate fields $K \subs E \subs L$:
$$
\HT^0(\Gal(L/E), \Z_B(n, L)) = 0,
$$
where $\HT^q(-,-)$ denotes Tate cohomology, see Remark \ref{ht0}.

By Hilbert's Theorem 90, we have that:
$$
\HT^1(\Gal(L/E), \Z_B(2, L)) = 0,
$$
for any intermediate field $E$, see Remark \ref{hthm90}.  Thus, by Corollary \ref{htwin}, we have that:
$$
\HT^q(\Gal(L/K), \Z_B(n, L)) = 0
$$
for {\em all} $q \in \Z$.  

Taking a limit over finite Galois extensions $L/K$, we see that:
$$
H^q(G_K, \Ze(n, K)) = 0
$$
for all $q>0$ and all finite extensions $K/k$.

Now, let $p$ be a prime different from the characteristic of $k$.  We have an exact triangle in the derived category of $G_K$-modules:
$$ 
\Ze(n, K) \ra \Ze(n, K) \ra \mu_p^{\otimes n}[n] \ra \Ze(n, K)[1].
$$
Thus, if $K/k$ is any finite extension, we have that:
$$
H^q(G_K, \mu_p^{\otimes n}) = 0
$$
for all $q > n$.  

Now, let $H$ be a $p$-Sylow subgroup of $G_k$. Let $S$ be the field corresponding to $H$.  The extension $S/k$ is a limit over all finite Galois extenions $K/k$ contained in $S$.  Taking a limit over these fields, we see that:
$$
H^{n+1}(H, \mu_p^{\otimes n}) = 0.
$$
But $H$ acts trivially on $\mu_p$, because $S$ contains the $p^{th}$ roots of unity.  Thus, we can identify the $n^{th}$ Tate twist $\mu_p^{\otimes n}$ with $\Z/p\Z$, obtaining:
$$
H^{n+1}(H, \Z/p\Z) = 0.
$$
But by Proposition 21 in $\sect$I.4.1 of Serre \cite{CG}, this implies that $cd_p(H) \leq n$, which implies that $cd_p(G_k) \leq n$. 

Finally, if $p$ is the characteristic of $k$, then $cd_p(G_k) \leq 1$, see Proposition 3 in $\sect$II.2.2 of Serre \cite{CG}.  Thus, the cohomological dimension of $k$ is at most $n$.    

Next, suppose that $cd(G_k) \leq n$.  If $K/k$ is any finite extension, then $G_K$ is isomorphic to a subgroup of finite index in $G_k$, so we have that 
$$
cd(G_K) \leq cd(G_k) \leq n, 
$$
see Proposition 22 in $\sect$I.4.1 of Serre \cite{CG}.  

For every prime $p$ not equal to the characteristic of $k$,
$$
H^2(G_K, \Ze(n, K))(p) = 0,
$$
by Lemma 1.1 of Saito \cite{S4}.  This can be easily deduced using the exact triangle:
$$
\Ze(n, K) \ra \Ze(n, K) \ra \mu_p^{\otimes n}[n] \ra \Ze(n, K)[1].
$$  

If the characteristic of $k$ is equal to $p>0$, then there is an exact triangle in the derived category of $G_K$-modules:
$$
\Ze(n, K) \ra \Ze(n, K) \ra v(n)_K \ra \Ze(n, K)[1],
$$
where $v(n)_K$ denotes the additive subgroup of the deRham complex 
$\Omega_K^{\bullet}$ generated by logarithmic differential $n$-forms, see Milne \cite{Mflat}.  But, we have that 
$$
[K:K^p] \leq [k:k^p] \leq p^{n-1}, 
$$
which implies that the dimension of $\Omega_K^1$ is at most $n-1$.  So, for completely elementary reasons, we have that $\Omega_K^n = 0$, and:
$$
H^2(G_K, \Ze(n, K))(p) = 0
$$
for $p$ equal to the characteristic of $k$.

Thus, for all finite extension $K/k$, we have that:
$$
H^2(G_K, \Ze(n, K)) = 0.
$$

Now, suppose that $L/K$ is a finite Galois extension.  We have an inflation-restriction exact sequence:
$$
0 \ra H^2(\Gal(L/K), \Z_B(n, L)) \ra H^2(G_K, \Ze(n, K)) \ra H^2(G_L, \Ze(n, L)),
$$
see Remark \ref{inflres}.  This implies that:
$$
H^2(\Gal(L/K), \Z_B(n, L)) = 0
$$
for all finite Galois extensions $L/K$ with $K/k$ finite.  Again, we can use Hilbert's Theorem 90 and Corollary \ref{htwin} to see that:
$$
\HT^q(\Gal(L/K), \Z_B(n, L)) = 0
$$
for all $q \in \Z$.  

In particular, setting $q=0$ gives the surjectivity of the norm map:
$$
N : K_2(L) \lra K_2(K).
$$
\end{proof}

\begin{remark} By setting $n=1$, one recovers Serre's proof in \cite{CG} of the Theorem stated in the introduction.
\end{remark}

\begin{remark}  The condition that $[k:k^p] \leq p^{n-1}$ if $k$ is of characteristic $p>0$ is only used to imply the triviality of $H^2(G_K, \Ze(n, K))(p)$.  One could replace this condition with the weaker  assumption that:
$$
H^2(G_K, \Ze(n, K))(p) = 0
$$
for all finite extensions $K/k$, as in Proposition 5 in $\sect$II.3.1 of Serre \cite{CG}.
\end{remark}

\begin{remark}  Using the notation of Kato introduced in \cite{K2} and \cite{Kato0}, one can rephrase the main result: if $k$ is any field, then $k$ is a $\mathfrak{B}_n$-field if and only if the dimension of $k$ is at most $n$. 
\end{remark}

\section{Applications of the main result}

\begin{app}  Suppose that $k$ is a global field, that is, a number field or the function field of a curve over a finite field.  Then it is well-known, see Serre \cite{CG}, that the cohomological dimension of $k$ is at most 2, except in the case that $k$ is a number field with a real place.  Furthermore, if $k$ is a function field, then we have that $[k:k^p] = p$.  Thus, for function fields and totally imaginary number fields, we recover a theorem of Bak \cite{Bak} on the surjectivity of the norm map on $K_2$.  

In the case that $k$ is a number field with at least one real embedding, then we have that $cd_p(G_k) \leq 2$ for all $p \neq 2$.  We also recover the theorem of Bak {\em loc. cit.} that for number fields with a real place, the norm map is surjective on $K_2$ except on the 2-primary part.  This can be explicitly seen using the description of the Milnor ring of a global field due to Bass-Tate \cite{BT}, and the fact that any local field is a $\mathfrak{B}_2$-field.
\end{app}

\begin{remark}  In the course of the proof of the main result, we saw that $k$ is a $\mathfrak{B}_n$-field if and only if for all finite extensions $K/k$, a certain motivic cohomology group vanishes; that is, that:
$$
H^2(G_K, \Ze(n, K)) = 0.
$$
For $n=1$, this motivic cohomology group is the Brauer group of $K$.  One should think of these cohomology groups for $n>1$ as the ``higher weight'' Brauer groups of $K$.  

One can rephrase the main result using this language: if $k$ is a field such that for all finite extensions $K/k$, the weight $n$ Brauer group of $K$ is trivial, then the cohomological dimension of $k$ is at most $n$.  The converse is true, after adding the additional assumption that either $[k:k^p] \leq p^{n-1}$, or $H^2(G_K, \Ze(n, K))(p) = 0$.
\end{remark}

\begin{app}  Let $X= \Spec R$, where $R$ be a number ring with field of fractions $K$, a totally imaginary number field.  As noted previously, $K$ is a $\mathfrak{B}_2$-field.  Let $\nu$ be a finite place of $K$, and let $K_{\nu}$ denote the completion of $K$ at $\nu$.  Then $K_{\nu}$ is a complete discrete valuation field with residue field $k_{\nu}$, a finite field.  Each $k_{\nu}$ is a $\mathfrak{B}_1$-field.  

So, $X$ is a scheme whose generic point is a $\mathfrak{B}_2$-field, and whose codimension one points are $\mathfrak{B}_1$-fields.  We have a long exact sequence of cohomology groups coming from localization theory:
\begin{eqnarray*}
\ldots \ra H^q(X, \Ze(n, X)) &\ra& H^q(G_K, \Ze(n, K)) \\ 
&\ra& \bigoplus_{\nu}H^q(G_{k_{\nu}}, \Ze(n-1, k_{\nu})) \ra \ldots.
\end{eqnarray*}
The main result implies that the motivic cohomology groups $H^q(X, \Ze(n, X))$ are trivial for $n>1$ and $q>1$.  One has the same result when $X$ is a regular curve over a finite field.
\end{app}

\begin{app}  In \cite{K2}, Kato proves that if $F$ is a complete discrete valuation field with residue field $k$, then $F$ is a $\mathfrak{B}_{n+1}$-field if and only if $k$ is a $\mathfrak{B}_n$-field.  One can use the main result to give a simple proof of this fact, using Kato's result \cite{Kato0} that the dimension of $F$ is equal to the dimension of $k$ plus one.
\end{app}

\begin{app} Classically, $\mathfrak{B}_1$-fields play an important role in local and global class field theory.  In the local theory, one uses the fact that this is a property of a finite fields, and of the maximal unramified extension of a local field.

In the global theory, one can use this property to understand the parallel significance of the cyclotomic extensions of a number field, and extensions of the base field of a function field.  The maximal such extension is both cases is a $\hat{\Z}$-extension, and the resulting field will be a $\mathfrak{B}_1$-field.
\end{app}

\begin{app}  As one might expect, $\mathfrak{B}_n$-fields for $n>1$ make an appearance in the cohomological approach to higher local class field theory, as in Kato \cite{K2}, Koya \cite{Ko}, Spiess \cite{Sp} and Saito \cite{S4}.  This concerns the study of the Galois group of $n$-local fields, that is, a complete discrete valuation field with residue field an $(n-1)$-local field, where a zero-local field is declared to be a finite field.  If $k$ is an $n$-local field, then $k$ is itself a $\mathfrak{B}_{n+1}$-field; the maximal unramified extension of $k$ in which all corresponding residue field extensions are unramified is the ``smallest'' field over $k$ which is a $\mathfrak{B}_n$-field.  

One expects that $\mathfrak{B}_2$-fields will play an important role in a cohomological approach to global class field theory for arithmetic surfaces, if such an approach exists.  One can see a manifestation of this phenomenon in the semi-global result of \cite{me}.  
\end{app}

\begin{app}  It may be possible to relate Kato's $\mathfrak{B}_n$ property to Serre's $C_n$ property \cite{CG}.  Recall that a field $k$ is called a $C_n$-field if every homogeneous polynomial of degree $d$ in more than $d^n$ variables has a nontrivial zero.  

Kato \cite{K2}*{$\sect$3} showed that every $C_2$-field is a $\mathfrak{B}_2$ using the reduced norm.  The converse is false, because it is well-known that field $\Q_2$ of 2-adic numbers is {\em not} a $C_2$-field, although its cohomological dimension is 2, which implies that it is a $\mathfrak{B}_2$-field.  

Kato conjectured in\cite{K2} that every $C_n$-field is a $\mathfrak{B}_n$-field; using the main result, we see that this is equivalent to Serre's conjecture in \cite{CG} that the cohomological dimension of a $C_n$-field is at most $n$.  
\end{app}

\appendix

\section{Motivic complexes} \label{motivic}

Let $X$ be a scheme.  In \cite{Li0}, Lichtenbaum predicted the existence of certain objects in the derived category of sheaves on the $\e$tale site over $X$, which have come to be known as motivic complexes. First, we present a list of properties (or ``Axioms'') that motivic complexes are expected to satisfy.  These properties have been modified slightly: we work with complexes that are acyclic in {\em positive} degrees.  The reader will find that the grading used here is shifted down by $n$ from the standard choice {\em loc. cit}.

\begin{enumerate}
\item If $n>0$, $\Ze(n, X)$ is acyclic outside the interval $[-n+1, 0]$.
\item Let $\al_*$ be the functor that assigns to every $\e$tale sheaf on $X$ the associated Zariski sheaf.  Then $\R^{1} \al_* \Ze(n, X)=0.$
\item Let $l$ be a positive integer, prime to all residue field characteristics of $X$.  Then there is an exact triangle:
$$ 
\Ze(n, X) \ra \Ze(n, X) \ra \mu_l^{\otimes n}[n] \ra \Ze(n, X)[1],
$$
where $\mu_l^{\otimes n}$ denotes the $n^{th}$ Tate twist of the $l^{th}$ roots of unity $\mu_l$ on the $\e$tale site over $X$.

If $p$ is a residue field characteristic of $X$, then there is an exact triangle in the derived category of $G_K$-modules:
$$
\Ze(n, X) \ra \Ze(n, X) \ra v(n) \ra \Ze(n, X)[1],
$$
where $v(n)$ denotes the additive subsheaf of the deRham complex 
$\Omega_X^{\bullet}$ generated by logarithmic differential $n$-forms, see Milne \cite{Mflat}.  
\item There are product maps:
$$
\Ze(n, X) \otimes^L \Ze(m, X) \lra \Ze(n+m, X).
$$
\item The sheaves $\Hc^i(\Ze(n, X))$ are isomorphic to $\gr_{\gamma}^n \K_{n-i}$, where $\K_*$ is the sheaf of algebraic $K$-groups on the $\e$tale site over $X$, and $\gr_{\gamma}$ is the gradation corresponding to Soul$\e$'s $\gamma$-filtration \cite{Soule}.
\item If $F$ is a field, then $H^0(G_F, \Ze(n, F))$ is canonically isomorphic to $K_n^M(F)$.  
\end{enumerate}

\begin{remark} For small weights, the motivic complex is concentrated in degree zero.  $\Ze(0, X)$ is the constant sheaf $\Z$; $\Ze(1, X)$ is the multiplicative group $\G_m$.  
\end{remark}

\begin{remark} \label{hthm90}  If $K$ is a field, then the second axiom for $X=\Spec K$ is an analogue of Noether's generalization of Hilbert's Theorem 90.  It is equivalent to stipulating that:
$$
H^1(G_K, \Ze(n, K)) = 0.
$$
\end{remark}

\begin{remark} \label{zb} If $K$ is a field, and $L/K$ is a Galois extension, then there is a natural way to produce a complex of $\Gal(L/K)$-modules from a complex of $G_K$-modules.  It is given by: 
$$
C^{\bullet} \mapsto \trg(G_L, C^{\bullet}).
$$
In the case that $C^{\bullet}$ is concentrated in degree zero, we recover the classical operation of taking $G_L$-invariants.

For the weight $n$ motivic complex over a field $K$, it is expected that the corresponding complex of $\Gal(L/K)$-modules will have much in common with the $\e$tale complex.  We denote
$$
\Z_B(n, L) = \trg(G_L, \Ze(n, K)),
$$
after Beilinson, who first envisioned motivic complexes for the Zariski site. 
\end{remark}

There is also a set of axioms for Beilinson's complex $\Z_B(n, X)$ of sheaves on the Zariski site of $X$:

\begin{enumerate}
\item If $n>0$, then $\Z_B(n, X)$ is acyclic outside of degrees $[-n+1, 0]$.
\item If $l$ is invertible on $X$, then 
$$
\Z_B(n, X) \otimes^L \Z/l\Z = \tn R\al_* \Z/l\Z(n)[n],
$$
where $\Z/l\Z(n)$ denotes the $n^{th}$ Tate twist of $\mu_l$.
\item $\gr^n_{\gamma} \K_j \simeq H^{n-j}(X, \Z_B(n, X))$ up to torsion, and probably up to ``small factorials''.  Here $\K_*$ denotes the sheaf of $K$-groups on the {\em Zariski} site of $X$.
\item If $X$ is smooth, then $H^0(X, \Z_B(n, X))$ is $\K_n^M(X)$.  
\end{enumerate}

In general, $\Z_B(n, X)$ should be given by $\tn \R\al_*\Ze(n, X)$.  In fact, Beilinson's axioms follow directly from a strengthened version of the axioms for the $\e$tale site, see Lichtenbaum \cite{Li0}.

\begin{remark}  \label{inflres}  Let $K$ be a field, and $L/K$ a finite Galois extension.  Then the cohomology groups for the Zariski motivic complex $\Z_B(n, L)$ and the $\e$tale motivic complex $\Ze(n, K)$ are related through an inflation-restriction exact sequence:
$$
0 \ra H^2(\Gal(L/K), \Z_B(n, K)) \ra H^2(G_K, \Ze(n, K)) \ra H^2(G_L, \Ze(n, L)).
$$
This is a simple consequence of Hilbert's Theorem 90 for $\Ze(n, K)$ and the fact that there is a quasi-isomorphism:
$$
\Z_B(n, L) = \trg(G_L, \Ze(n, K)).
$$
\end{remark}

\begin{remark} In 1987, Lichtenbaum \cite{Li} constructed a candidate for the weight two motivic complex, using techniques from $K$-theory.  He showed it satisfies most of the Axioms, and all the Axioms if $X = \textup{Spec } K$.  

Bloch was the first to construct a candidate for the weight $n$ motivic complex, using higher Chow groups, see \cite{Bl}.  Using another approach, Voevodsky, Suslin and Friedlander \cite{VSF} have constructed a triangulated category that plays the role of the derived category of mixed motives.  In the process, they have shown \cite{V} that Bloch's higher Chow group complex satisfies all of Lichtenbaum's axioms, except perhaps the first.  This is axiom is equivalent to the Beilinson-Soul$\e$ Vanishing Conjecture \cite{Soule}.  It is known that Soul$\e$'s $\gamma$-filtration has finite length for $X = \textup{Spec } K$, that is, that Bloch's complex is bounded below.  For the purpose of our result on cohomological dimension, this is sufficient.  
\end{remark}

\section{Tate cohomology} \label{tate}

When $G$ is a finite group, we use the standard two-sided resolution $P_{\bullet}$ of $\Z$ as a $G$-module to define Tate cohomology groups, see Neukirch-Schmidt-Wingberg \cite{NSW}.   The following definition was made by Koya \cite{Ko}:

\begin{definition} Let $G$ be a finite group, $A^{\bullet}$ a bounded complex of $G$-modules, and $P_{\bullet}$ the standard two-sided resolution of $\Z$ as a $G$-module.  The Tate cohomology group $\HT^q(G, A^{\bullet})$ for $q \in \Z$ is defined to be the $q^{th}$ total homology of the double complex: 
$$
\Hom_G(P_{\bullet}, A^{\bullet}).
$$
\end{definition}

\begin{remark} \label{norms}  In this setting, one can introduce norm maps on hyper-cohomology groups.  Let $N = \sum_{\sig \in G} \sig$ be the norm element of the group ring $\Z[G]$.  If $A^{\bullet}$ is a complex of $G$-modules, then multiplication by $N$ yields a complex $NA^{\bullet}$ with the property that 
$$
\Hc^q(NA^{\bullet}) = N \Hc^q(A^{\bullet}).
$$
The latter group can be considered as a subgroup of $H^0(G, \Hc^q(G, A^{\bullet}))$, which maps via an edge morphism to $H^q(G, A^{\bullet})$, using the standard spectral sequence to compute hyper-cohomology, see Weibel \cite{W}.  The norm map on hyper-cohomology groups is defined the be the composite:
$$
N : \Hc^q(A^{\bullet}) \lra H^q(G, A^{\bullet}).
$$

In the case where $G$ is  $\Gal(L/K)$ and $A^{\bullet}$ is $\Z_B(n, L)$, the weight $n$ motivic complex on the Zariski site, we can recover the norm maps on Milnor $K$-theory defined by Bass-Tate \cite{BT} by setting $q=0$.  These definitions necessarily agree, because the definition of the norm through hyper-cohomology shares the same functoriality properties as in Kato \cite{K2}*{$\sect$1}, which characterize the norm map uniquely.
\end{remark}

\begin{remark} \label{ht0}  For the purposes of Tate cohomology, the assumption that $A^{\bullet}$ is a bounded complex, that is, that $A^{\bullet}$ has only finitely many non-zero entries, is necessary for the convergence of the spectral sequence associated to the double complex.  
If $A^{\bullet}$ is a complex that is zero in degree larger than $n$, then 
$$
\HT^q(G, A^{\bullet}) \simeq \Hc^q\big(\RG(G, A^{\bullet})\big) = H^q(G, A^{\bullet})
$$
for all $q > n$.  Furthermore, we have that 
$$
\HT^n(G, A^{\bullet}) = \Hc^n\big(\RG(G, A^{\bullet})\big)/N\Hc^n\big(A^{\bullet}\big), 
$$
see Theorem 1.2 of Koya \cite{Ko}.
\end{remark}

We have the following fundamental observation:

\begin{prop}[Dimension-shifting for complexes] \label{cshift} Let $G$ be a finite 
group and $A^{\bullet}$ a bounded 
complex of $G$-modules.  Then there exists a $G$-module $A'$ such that for 
any subgroup $H$ of $G$, 
$$
\HT^q(H, A^{\bullet}) \simeq \HT^q(H, A')
$$
for {\em all} $q \in \Z$.
\end{prop}

See \cite{weil} for a proof.

\begin{cor} \cite{Ko1}*{Prop. 4}  \label{htwin} Let $G$ be a finite group, and $A^{\bullet}$ a bounded complex of $G$-modules.  Suppose that for every $p$-Sylow subgroup $H$ of $G$, there is an integer $i_H$ such that
$$
\HT^{i_H}(H, A^{\bullet}) = \HT^{i_H+1}(H, A^{\bullet}) = 0.
$$
Then $\HT^q(G, A^{\bullet}) = 0$ for all $q \in \Z$.
\end{cor}

This follows immediately from the Proposition, and the corresponding classical fact for the Tate cohomology of a module, see Artin-Tate \cite{AT}.

\end{document}